\documentclass[12pt]{amsart}

\usepackage{amssymb}

\usepackage{enumerate}

\usepackage{graphicx}

\makeatletter
\@namedef{subjclassname@2010}{%
  \textup{2010} Mathematics Subject Classification}
\makeatother

\newtheorem{thm}{Theorem}[section]

\newtheorem{prop}[thm]{Proposition}

\theoremstyle{definition}
\newtheorem{defin}[thm]{Definition}
\newtheorem{rem}[thm]{Remark}
\newtheorem{exa}[thm]{Example}

\newtheorem*{prf}{Proof}
\numberwithin{equation}{section}

\begin{document}

%%%%% To ease editing, for IMPAN journals add:

\baselineskip=17pt

%%%%%%%%%%%

%% In the running head, replace first names by initials 
%% and give an abbreviation of the title.

\title[Sample paper]{On Cubic Residue Matrices}

\author[R. Wood]{Ryan Wood}
\address{Department of Mathematics and Statistics\\ Northern Arizona University\\
Flagstaff, Arizona 86011}
\email{rpw54@nau.edu}

\author[J. Rushall]{Jeff Rushall}
\address{Department of Mathematics and Statistics\\ Northern Arizona University\\
Flagstaff, Arizona 86011}
\email{jeffrey.rushall@nau.edu}

\author[P. Gonzalez]{Pauline Gonzalez}
\address{Department of Mathematics and Statistics\\ Northern Arizona University\\
Flagstaff, Arizona 86011}
\email{pg347@nau.edu}

\date{}

\begin{abstract}
The use of quadratic residues to construct matrices with specific determinant values is a familiar problem with connections to many areas of mathematics and statistics. Our research has focused on using \emph{cubic residues} to construct matrices with interesting and predictable determinants. 
%A template for articles in IMPAN journals in the \texttt{amsart} style. Using \texttt{pdflatex} is strongly preferred.
\end{abstract}



\maketitle

\section{Introduction}

%Given any odd prime number $p$, consider the quadratic congruence equation ${x}^2 \equiv a $(mod$ \ p$). If there exists a value of $x$ that satisfies this congruence equation, then $a$ is a quadratic residue modulo $p$, whereas if no solution exists, $a$ is a quadratic nonresidue modulo $p$. It is convenient to reformulate these concepts with the help of the \emph{Legendre symbol modulo} $p$ (also known as the \emph{quadratic residue symbol modulo} $p$), which is defined as:
%\[ \bigg( \frac{a}{p} \bigg)= \begin{cases} 
%      1 & \text{if}\; \textit{a}\; \text{is a nonzero quadratic residue modulo}\; \textit{p} \\
%      -1 & \text{if}\; \textit{a}\; \text{is a quadratic nonresidue modulo}\; \textit{p} \\
%      0 & \text{if}\; a \equiv 0 (\text{mod} \ p) 
%   \end{cases}
%\]

%\vspace{3mm}

In 1933, R. Paley [3] used Legendre symbols in an ingenious way to construct two new infinite classes of Hadamard matrices, once each for primes of the form $4k+1$ and $4k+3$. For instance, the Paley construction associated with $4k+3$ primes uses matrices whose $i^{th}$ row, $j^{th}$ column entry is defined by the quadratic residue symbol $\big( \frac{i-j}{p} \big)$.

R. Chapman [1] has studied matrices generated by various types of Legendre symbols. For instance, consider the matrix $M=[m_{ij}]$ whose $i^{th}$ row, $j^{th}$ column entry is defined by $m_{ij}$= $\big( \frac{j-i}{p} \big)$. Chapman observed that if $p$ is a small prime of the form $4k+3$, then the matrix $M$ of order $\frac{p+1}{2}$ defined as above satisfied $det(M)=1$. Chapman subsequently verified this for all such primes $p<1000$, but was unable to prove this held for all $4k+3$ primes, a claim which soon became known as Chapman's Evil Determinant Problem. His conjecture was verified by M. Vsemirnov in 2012 [4]; Vsemirnov later proved a related conjecture involving matrix determinants based on primes of the form $4k+1$ [5]. 

Both Paley and Chapman used linear relationships in terms of $i$ and $j$ in their Legendre symbol formulas to generate matrices with interesting and/or useful determinants. In this paper, we attempt to model the approaches of Paley and Chapman to build matrices with predictable determinants, but we generate matrix entries using cubic residue symbol values rather than quadratic residue symbol values.

%You can use this file as a template when submitting your paper to one of IMPAN journals (except Dissertationes Mathematicae and Banach Center Publications, for which style files exist).

%The format of this file is \textbf{not} the exact final printed format (the latter is scaled down, and line breaks will most often be different), but it is convenient for editing purposes.

\section{Background Information}
For the remainder of this paper, all variables are assumed to be integers, and all moduli are assumed to be integer primes.

\begin{defin} If there exists an integer solution to the equation ${x}^3 \equiv a $(mod$ \ m$), we say that $a$ is a \emph{cubic residue modulo} $m$. If no solution exists, $a$ is a \emph{cubic nonresidue modulo} $m$.
\end{defin}

%\begin{exa}
Since ${4}^3 = 64 \equiv 12 $(mod$ \ 13$), $12$ is a cubic residue modulo $13$. Via brute force, we find that ${x}^3 \equiv 2 $(mod$ \ 7$) has no solutions, and so $2$ is a cubic nonresidue modulo $7$.
%\end{exa}

\begin{defin}
Given any odd prime $p$, we define the \emph{cubic residue symbol modulo} $p$ thusly:

\[ \bigg[ \frac{a}{p} \bigg]= \begin{cases} 
      1 & \text{if}\; \textit{a}\; \text{is a nonzero cubic residue modulo}\; \textit{p} \\
      -1 & \text{if}\; \textit{a}\; \text{is a cubic nonresidue modulo}\; \textit{p} \\
        0 & \text{if}\; a \equiv 0 (\text{mod} \ p)

   \end{cases}
\]
\end{defin}

%\begin{exa}

For instance, since $12$ is a cubic residue modulo $13$, $\big[ \frac{12}{13} \big]= 1$, and since $2$ is a cubic nonresidue modulo $7$, $\big[ \frac{2}{7} \big]= -1$.

%\end{exa}

The cubic residue symbol has several useful and well-known properties which are summarized in the following result. 

\begin{prop}
Given any odd prime $p$ and any integers $a$ and $b$, then:
\begin{itemize}

\item if $a \equiv b $(mod$ \ p$), then $\big[ \frac{a}{p} \big]= \big[ \frac{b}{p} \big]$

\item if $a \not\equiv 0$(mod$ \ p$), then $\big[ \frac{a^3}{p} \big]=1$

\item $\big[ \frac{a}{p} \big]= \big[ \frac{-a}{p} \big]$

\end{itemize}
\end{prop}

The behavior of the values of a cubic residue symbol modulo $p$ depends on whether the prime is of the form $3k+1$ or $3k+2$. If $p=3k+1$, then the distribution of the corresponding cubic residue symbol values $1$ and $-1$ is essentially random for small values of $a$. But the number of cubic residues is relatively easy to determine.

\begin{prop}

If $p$ is a $3k+1$ prime, then among the $p-1$ nonzero congruence classes there are exactly $\frac{p-1}{3}$ cubic residues and $\frac{2(p-1)}{3}$ cubic nonresidues. 

\end{prop}

On the other hand, if $p=3k+2$, we have the following result (see, for instance, [2]).
\begin{prop}
If $p$ is a $3k+2$ prime, then $\big[ \frac{a}{p} \big]= 1$ for every nonzero $a$. 
\end{prop}

\begin{exa}
The cubic residues modulo $7$ are $1$ and $6$; the nonresidues are $2,3,4$ and $5$. Here is the order $3$ matrix $M=$ [$m_{ij}$], whose $i^{th}$ row, $j^{th}$ column entry is defined by $m_{ij}$= $\big[ \frac{j-i}{7} \big]$:

\[
\begin{bmatrix}
    [\frac{1-1}{7}] & [\frac{2-1}{7}] & [\frac{3-1}{7}] \\
    [\frac{1-2}{7}] & [\frac{2-2}{7}] & [\frac{3-2}{7}] \\
    [\frac{1-3}{7}] & [\frac{2-3}{7}] & [\frac{3-3}{7}] \\
\end{bmatrix}
=
\begin{bmatrix}
    [\frac{0}{7}] & [\frac{1}{7}] & [\frac{2}{7}] \\
    [\frac{-1}{7}] & [\frac{0}{7}] & [\frac{1}{7}] \\
    [\frac{-2}{7}] & [\frac{-1}{7}] & [\frac{0}{7}] \\
\end{bmatrix}
=
\begin{bmatrix}
    0 & 1 & -1 \\
    1 & 0 & 1 \\
    -1 & 1 & 0 \\
\end{bmatrix}
\]
\end{exa}

\section{New Results}

Our plan is to construct matrices whose entries are specific cubic residue symbol values, compute their determinants, organize the data and look for patterns. We have several options, which include:
\begin{itemize}
\item choosing either the $3k+1$ or $3k+2$ prime category;

\item choosing the prime $p$ on which the cubic residue symbol is based; 

\item choosing the order $n$ of the matrix;

\item choosing the $i,j$ "formula" within the cubic residue symbol;

\item choosing how to organize the resulting determinant values.
\end{itemize}

This process was not unlike looking for a needle in a field of haystacks without knowing if a needle even existed. We chose to first focus on matrices generated by cubic residue symbols modulo primes of the form $3k+1$, in part because the robustness of these particular cubic residue values closely parallel those of quadratic residue values. 

In Figure 1 we see a typical table of matrix determinants, color-coded for ease of analysis. This determinant table results from using the cubic residue symbol formula $m_{ij} = [\frac{j+i+c}{7}]$. The vertical axis describes the matrix order $n$; the horizontal axis corresponds to $c-$values; the color coding has blue denoting zero determinants, orange denoting negative determinants, and green denoting positive determinants. \\

\begin{figure}[h]
%\begin{minipage}{0.5\textwidth}\centering
\includegraphics[width=\linewidth]{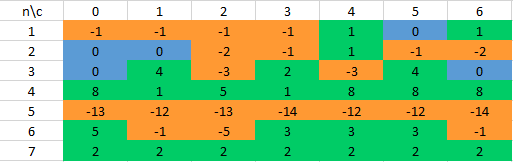}
\caption{The determinant table generated by $m_{ij} = [\frac{j+i+c}{7}]$.}
\end{figure}
%\end{minipage}
%\begin{minipage}{0.4\textwidth}\raggedleft

The pseudo-random distribution of determinant values and $+/-$ \, signs is fairly representative of the determinant tables that result from using $3k+1$ primes. The lack of a perceivable pattern in this and the many other $3k+1$ tables we generated led us to consider $3k+2$ primes.
%\end{minipage}
%\noindent
%\\

The determinant table that results from using the $3k+2$ prime $p=11$ and $m_{ij} = [\frac{j-i+c}{11}]$ appears in Figure 2. We color-coded determinant values to highlight the surprising visual structure of this table. The pattern appearing in this table repeats horizontally and features constant zero determinant values for all matrices of order $n>p$ and for every prime $p$ of the form $p=3k+2$, as seen in Figure 3. The table structure allows us to formulate and prove several new results about cubic residue matrix determinants.

\begin{figure}[h]
\begin{center}
\includegraphics[width=\linewidth]{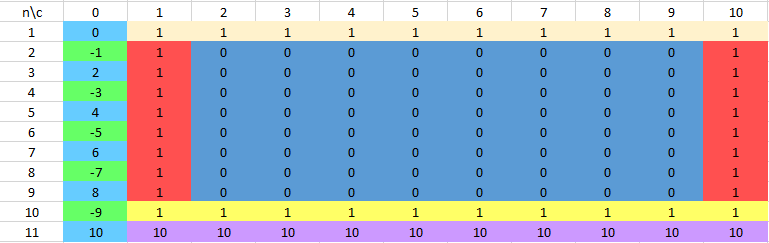}\\
\end{center}
\caption{The determinant table generated by $m_{ij} = [\frac{j-i+c}{11}]$.}
\end{figure}

\begin{figure}[h]
\begin{center}
\includegraphics[width=\linewidth]{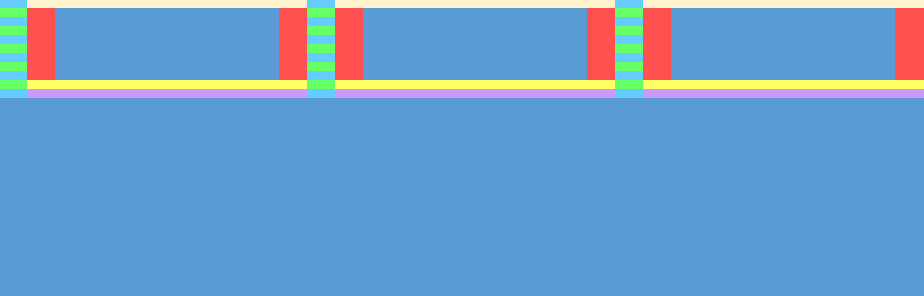}\\
\end{center}
\caption{The extended table generated by $m_{ij} = [\frac{j-i+c}{11}]$.}
\end{figure}

\begin{thm}
Given any prime $p$ of the form $3k+2$, let the matrix $M_n$ of order $n$ be defined by $M_{n}=[m_{ij}]$ with $m_{ij}$ = $\big[ \frac{j-i+c}{p} \big]$. Then: $det(M_n)=(-1)^{n-1}(n-1)$ if $1 \leq n \leq p$ and $c=0$. 
\end{thm}
\begin{rem}
Note that this result describes the bicolor pattern in the initial column of the table in Figure 2.
\end{rem}
\begin{prf}
For each value of $n$, the matrix $M_n$ generated by these cubic residues is simply the adjacency matrix of a complete graph on $n$ vertices, and the result follows. 
\openbox
\end{prf}

%\vspace{8mm}

\begin{thm}
Given any prime $p$ of the form $3k+2$, let the matrix $M_n$ of order $n$ be defined by $M_{n}=[m_{ij}]$ with $m_{ij}$ = $\big[ \frac{j-i+c}{p} \big]$. Then $det(M_n)=1$ if $1 < n \leq p-1$ and $c=\pm 1$. 
\end{thm}
\begin{rem}
This result describes the values that appear in both the second column and the last column of the table in Figure 2.
\end{rem}
\begin{prf}
Let $c=1$. Clearly, $det(M_2)=1$. Now assume that $n>2$ and let $R_1,R_2,\ldots , R_n$ denote the rows of $M_n$. Note that $M_n$ has $1$'s for every entry except for $0$ entries on the subdiagonal. If $n=3$, $M_3$ can be made upper triangular - with all nonzero entries of $1$ - by replacing $R_1$ with $R_1+R_2-R_3$, and so $det(M_3)=1$.

We now proceed inductively. For any $n>3$, the following row operations transform $M_n$ into an upper triangular matrix with all nonzero entries of $1$: 
\vspace{2mm}
\begin{center}
$R_1 = \displaystyle\sum_{i=1}^{n-1} R_i - R_n$, \\
$R_2 = \displaystyle\sum_{i=2}^{n-1} R_i - R_n$, \\
\end{center}
and similarly, for $3 \leq j \leq n-2$, \\
\begin{center}
$R_j = \displaystyle\sum_{i=j}^{n-1} R_i - R_n$. 
\end{center}
\vspace{2mm}
This process always yields an upper triangular matrix with all nonzero entries of $1$, and the result then follows. The case for $c=-1$ is similar and left to the reader.
\openbox
\end{prf}
%\vspace{8mm}

\begin{thm}
Given any prime $p$ of the form $3k+2$, let the matrix $M_n$ of order $n$ be defined by $M_{n}=[m_{ij}]$ with $m_{ij}$ = $\big[ \frac{j-i+c}{p} \big]$. Then $det(M_n)=p-1$ if $n=p$ and if $1 \leq c \leq p-1$. 
\end{thm}
\begin{rem}
This result describes the values that are in the bottom row of the table in Figure 2.
\end{rem}
\begin{prf}
When $c=0$ we obtain our result from Theorem 3.1; call this matrix $M_{p}^{*}$. Every value of $c$, $1 \leq c \leq p-1$, results in a matrix whose columns are an order-preserving translation of the columns of $M_{p}^{*}$. This rearrangement is achieved using $c(p-c)$ column transpositions, and since $p$ is an odd prime, exactly one of $c$ and $p-c$ is even. Hence, the determinant value remains unchanged.
\openbox
\end{prf}
%\vspace{8mm}

\begin{thm}
Given any prime $p$ of the form $3k+2$, let the matrix $M_n$ of order $n$ be defined by $M_{n}=[m_{ij}]$ with $m_{ij}$ = $\big[ \frac{j-i+c}{p} \big]$. Then $det(M_n)=0$ if $2 \leq n \leq p-2$ and $2 \leq c \leq p-2$. 
\end{thm}
\begin{rem}
This result describes the large rectangular block of $0$'s in the center of the table in Figure 2.
\end{rem}
\begin{prf}
For each value of $c$, $2 \leq c \leq p-2$, the resulting matrix $M_n$ features at least two columns of constant $1$'s for entries.
\openbox
\end{prf}
%\vspace{8mm}

\begin{thm}
Given any prime $p$ of the form $3k+2$, let the matrix $M_n$ of order $n$ be defined by $M_{n}=[m_{ij}]$ with $m_{ij}$ = $\big[ \frac{j-i+c}{p} \big]$. Then $det(M_n)=1$ if $n=p-1$ and $1 \leq c \leq p-1$. 
\end{thm}
\begin{rem}
This result describes the row of $1$'s near the bottom of the table in Figure 2.
\end{rem}

\begin{exa}
To understand how we came to this result, consider the matrix $M_{10}=[m_{ij}]$ with $m_{ij}$ = $\big[ \frac{j-i+4}{11} \big]$:
\[
\setcounter{MaxMatrixCols}{20}
\begin{bmatrix}
    1 & 1 & 1 & 1 & 1 & 1 & 1 & 0 & 1 & 1 \\
    1 & 1 & 1 & 1 & 1 & 1 & 1 & 1 & 0 & 1 \\
    1 & 1 & 1 & 1 & 1 & 1 & 1 & 1 & 1 & 0 \\
    1 & 1 & 1 & 1 & 1 & 1 & 1 & 1 & 1 & 1 \\
    0 & 1 & 1 & 1 & 1 & 1 & 1 & 1 & 1 & 1 \\
    1 & 0 & 1 & 1 & 1 & 1 & 1 & 1 & 1 & 1 \\
    1 & 1 & 0 & 1 & 1 & 1 & 1 & 1 & 1 & 1 \\
    1 & 1 & 1 & 0 & 1 & 1 & 1 & 1 & 1 & 1 \\
    1 & 1 & 1 & 1 & 0 & 1 & 1 & 1 & 1 & 1 \\
    1 & 1 & 1 & 1 & 1 & 0 & 1 & 1 & 1 & 1
\end{bmatrix}
\]. \\
Performing the column addition $C_{n}-C_{7} \rightarrow C_{n}$ for $n \neq 7$ we obtain the matrix:
\[
\setcounter{MaxMatrixCols}{20}
\begin{bmatrix}
    0 & 0 & 0 & 0 & 0 & 0 & 1 & -1 & 0 & 0 \\
    0 & 0 & 0 & 0 & 0 & 0 & 1 & 0 & -1 & 0 \\
    0 & 0 & 0 & 0 & 0 & 0 & 1 & 0 & 0 & -1 \\
    0 & 0 & 0 & 0 & 0 & 0 & 1 & 0 & 0 & 0 \\
    -1 & 0 & 0 & 0 & 0 & 0 & 1 & 0 & 0 & 0 \\
    0 & -1 & 0 & 0 & 0 & 0 & 1 & 0 & 0 & 0 \\
    0 & 0 & -1 & 0 & 0 & 0 & 1 & 0 & 0 & 0 \\
    0 & 0 & 0 & -1 & 0 & 0 & 1 & 0 & 0 & 0 \\
    0 & 0 & 0 & 0 & -1 & 0 & 1 & 0 & 0 & 0 \\
    0 & 0 & 0 & 0 & 0 & -1 & 1 & 0 & 0 & 0
\end{bmatrix}
\]

We now transpose column 7 with column 6, then with column 5, then with column 4, etc. until column 1 is a column of 1s. Replacing $C_{1}$ with the column sum $\sum_{k=1}^{n} C_{k}$ yields this matrix:

\[
\begin{bmatrix}
    0 & 0 & 0 & 0 & 0 & 0 & 0 & -1 & 0 & 0 \\
    0 & 0 & 0 & 0 & 0 & 0 & 0 & 0 & -1 & 0 \\
    0 & 0 & 0 & 0 & 0 & 0 & 0 & 0 & 0 & -1 \\
    1 & 0 & 0 & 0 & 0 & 0 & 0 & 0 & 0 & 0 \\
    0 & -1 & 0 & 0 & 0 & 0 & 0 & 0 & 0 & 0 \\
    0 & 0 & -1 & 0 & 0 & 0 & 0 & 0 & 0 & 0 \\
    0 & 0 & 0 & -1 & 0 & 0 & 0 & 0 & 0 & 0 \\
    0 & 0 & 0 & 0 & -1 & 0 & 0 & 0 & 0 & 0 \\
    0 & 0 & 0 & 0 & 0 & -1 & 0 & 0 & 0 & 0 \\
    0 & 0 & 0 & 0 & 0 & 0 & -1 & 0 & 0 & 0
\end{bmatrix}
\]
Appropriate column transpositions yield the following triangular matrix: 
\[
\begin{bmatrix}
    -1 & 0 & 0 & 0 & 0 & 0 & 0 & 0 & 0 & 0 \\
    0 & -1 & 0 & 0 & 0 & 0 & 0 & 0 & 0 & 0 \\
    0 & 0 & -1 & 0 & 0 & 0 & 0 & 0 & 0 & 0 \\
    0 & 0 & 0 & 1 & 0 & 0 & 0 & 0 & 0 & 0 \\
    0 & 0 & 0 & 0 & -1 & 0 & 0 & 0 & 0 & 0 \\
    0 & 0 & 0 & 0 & 0 & -1 & 0 & 0 & 0 & 0 \\
    0 & 0 & 0 & 0 & 0 & 0 & -1 & 0 & 0 & 0 \\
    0 & 0 & 0 & 0 & 0 & 0 & 0 & -1 & 0 & 0 \\
    0 & 0 & 0 & 0 & 0 & 0 & 0 & 0 & -1 & 0 \\
    0 & 0 & 0 & 0 & 0 & 0 & 0 & 0 & 0 & -1
\end{bmatrix}
\]
We see that the number of column transpositions is $1*6+3*7=27$, an odd integer. Thus, det$(M_{10}) = 1$.
\end{exa}

\begin{prf}
When $c=1$ we obtain our result from Theorem 3.3. Now consider $c$ such that $2 \leq c \leq p-1$. In this case, the matrix defined by $M_{n}=[m_{ij}]$ with $m_{ij}$ = $\big[ \frac{j-i+c}{p} \big]$ will be row equivalent to the following matrix:
\[
\setcounter{MaxMatrixCols}{20}
\begin{bmatrix}
    0 & 0 & \hdotsfor{2} & 0 & 1 & -1 & 0 & \dots & 0 & 0 \\
    \vdots & \vdots & \ddots & \ddots & \vdots & \vdots & 0 & \ddots & \ddots & \vdots & \vdots \\
    \vdots & \vdots & \ddots & \ddots & \vdots & \vdots & \vdots & \ddots & \ddots & 0 & \vdots \\
    \vdots & \vdots & \ddots & \ddots & \vdots & \vdots & \vdots & \ddots & \ddots & -1 & 0 \\
    0 & \vdots & \ddots & \ddots & \vdots & \vdots & \vdots & \ddots & \ddots & 0 & -1\\
    -1 & 0 & \ddots & \ddots  & \vdots & \vdots & \vdots & \ddots & \ddots & \vdots & 0\\
    0 & -1 & \ddots & \ddots  & \vdots & \vdots & \vdots & \ddots & \ddots & \vdots & \vdots \\
    \vdots & 0 & \ddots & \ddots  & \vdots & \vdots & \vdots & \ddots & \ddots & \vdots & \vdots \\
    \vdots & \vdots & \ddots & \ddots  & 0 & \vdots & \vdots & \ddots & \ddots & \vdots & \vdots \\
    0 & 0 & \dots & 0  & -1 & 1 & 0 & \hdotsfor{2} & 0 & 0
\end{bmatrix}
\]
where the $(p-c)^{th}$ column is the column of 1's. Using $p-c-1$ appropriate column transpositions, we move the $c^{th}$ column to the $(p-1)^{th}$ column. Then, after $(p-c)(c-1)$ more appropriate column transpositions, and recognizing that the column of 1's is equivalent to a column of 0's with a 1 in the $c^{th}$ row entry, we are left with triangular matrix whose main diagonal is comprised of $p-1$ $(-1)$'s. The result follows.
\openbox
\end{prf}
\begin{rem}
Given any prime $p$ of the form $3k+2$, the matrix $M_1$ of order $1$ defined by $M_{n}=[m_{ij}]$ with $m_{ij}$ = $\big[ \frac{j-i+c}{p} \big]$, for all $c > 0$ we have $det(M_1)=1$ because the only entry of the matrix is $1$.
\end{rem}

\begin{thm}
Given any prime $p$ of the form $3k+2$, for each order $2 \leq n \leq p-2$, define the matrix $M_n=[m_{ij}]$ via $m_{ij}$ = $\big[ \frac{j-i+1}{p} \big]$, and define the matrix $M^{'}_n=[m^{'}_{ij}]$ via $m^{'}_{ij}$ = $\big[ \frac{(j-i)^3+1}{p} \big]$. Then for all values of $i$ and $j$, $m_{ij} = m^{'}_{ij}$, and consequently $det(M_n)=det(M^{'}_n)$.
\end{thm}

\begin{prf}
For every prime $p$ of the form $3k+2$, the function $x \mapsto x^3$ is an automorphism of $\mathbb{Z}_{p}$ with $-1$ and $0$ fixed. Consequently, the cubic residue symbol values for $x+1$ and $x^3+1$ are identical modulo $p$. Thus, for all values of $i$ and $j$, 
\begin{center}
$\big[ \frac{j-i+1}{p} \big] = \big[ \frac{(j-i)^3+1}{p} \big]$, 
\end{center}
and the result follows.
\openbox
\end{prf}
%\vspace{5mm}

\begin{thm}
Given any prime $p$ of the form $3k+2$, construct $M=[m_{ij}]$, with $m_{ij}$ = $\big[ \frac{(j-i)^{2t}+c}{p} \big]$ of order $>1$. Let $r_p$ be a primitive root modulo $p$. Then $det(M)=0$ if both $p=12k+5$ and $c = r_{p}^{2n-1}$, or if both $p=12k+11$ and $c = r_{p}^{2n}$.
\end{thm}

% \begin{figure}
% \begin{center}
% \includegraphics[scale=1]{rectangleblue.PNG}
% \end{center}
% \caption{A monochrome table of determinants}
% \end{figure}

% %\vspace{4mm}
% This monochrome rectangle depicts the second result above, which describes the large block of nonsingular matrices that both of these particular cubic residue symbol formulas generate.\\

\begin{prf}

Let $p=12k+5$. Now suppose that $M$ contains a zero entry; that is, suppose that $(j-i)^{2t}+r_{p}^{2n-1} \equiv 0 (\text{mod} \ p)$ where $r_p$ is a primitive root modulo $p$ and $n,t \in \mathbb{Z}^{+}$. Then $(j-i)^{2t} \equiv -r_{p}^{2n-1} (\text{mod} \ p)$. Because $2n-1$ is odd, $(j-i)^{2t} \equiv (-r)_{p}^{2n-1} (\text{mod} \ p)$. Since $p \equiv 1 (\text{mod} \ 4)$, $-r_{p}$ is a primitive root modulo $p$. Then $(j-i)^{2t} \equiv ((-r)_{p}^{s})^{2t} \equiv (-r)^{2st} (\text{mod} \ p)$. But then $(-r)_{p}^{2n+1} \equiv (-r)^{2st} (\text{mod} \ p)$, a contradiction.\\
\indent Now let $p=12k+11$, and suppose that $(j-i)^{2u}+r_{p}^{2m} \equiv 0 (\text{mod} \ p)$ where $r_p$ is a primitive root modulo $p$ and $m,u \in \mathbb{Z}^{+}$ . Then $(j-i)^{2u} \equiv -r_{p}^{2m} (\text{mod} \ p)$. Using properties of Legendre symbols, and noting that $p \equiv 3 (\text{mod} \ 4)$, we see that 
\vspace{2mm}
\begin{center}
$1 = \Big( \frac{(j-i)^{2u}}{p} \Big) = \Big( \frac{-r_{p}^{2m}}{p} \Big) = \Big( \frac{r_{p}^{2m}}{p} \Big)\Big( \frac{-1}{p} \Big) = -\Big( \frac{r_{p}^{2m}}{p} \Big) = -1$, 
\end{center}
\vspace{2mm}
a contradiction, and so every entry of $M$ is nonzero. \\
\indent Given that $p \equiv 2 (\text{mod} \ 3)$, it must follow in both cases that each entry of $M$ will be $1$, and thus $det(M)=0$. \openbox
\end{prf}

\section{Future work}
We hope to generalize our work to the ring of Eisenstein integers, which is defined as $\mathbb{Z}[\omega]=\lbrace a+b\omega | a,b \in \mathbb{Z} \rbrace$, where $\omega$ is a primitive cubic root of unity. This generalization will also involve the use of a more sophisticated version of the cubic residue symbol than the one appearing in this paper.

%Preferably, figures should be prepared as pdf, jpg or eps files.
%All figures will be printed black and white; colours will only appear in the online version.

%Avoid very thin lines, and check whether all fonts used are embedded.

%Remember that sometimes figures have to be scaled, and then the lettering is scaled too; therefore, very small lettering should be avoided.

%\begin{figure}[h]
%\caption{A figure caption}
%\end{figure}

%\subsection*{Acknowledgements}
%This research was partly supported by NSF (grant no. XXXX).


\begin{thebibliography}{HD}

%% Use the widest label as the parameter.
%% Reference items can be numbered or have labels of your choice, as below.

%% In IMPAN journals, only the title is italicized; boldface is not used.
%% Our software will add links to many articles; for this, enclosing volume numbers in { } is helpful
%% Do not give the issue number unless the issues are paginated separately.



%%%%%%% To ease editing, use normal size:

\normalsize
\baselineskip=17pt

%%%%%%%%%%%%%%%
%M. Vsemirnov, \emph{On the Evaluation of R. Chapman's Evil Determinant}, Linear Algebra and its Applications, 436 (2012), 4101-4106.


\bibitem[1]{Chap}  R. Chapman,
My evil determinant problem (December, 2012),

\bibitem[2]{IreRos}  K. Ireland and M. Rosen,
\emph{Cubic and Biquadratic Reciprocity},
A Classical Introduction to Modern Number Theory {2} (1990), 108--137.

\bibitem[3]{Pal}  R. Paley,
\emph{On orthogonal matrices},
Journal of Mathematics and Physics {12} (1933), 311--320.

\bibitem[4]{Vsemirnov} M. Vsemirnov,
\emph{On the evaluation of R. Chapman's ``evil determinant''},  
in: Linear Algebra and its Applications (June, 2012),

\bibitem[5]{Vsemirnov2} M. Vsemirnov,
\emph{On R. Chapman's "evil determinant": case p=1 (mod 4)},
arXiv:1108.4031v2.









%\bibitem[K]{Kow}  J. Kowalski,
%\emph{A very interesting paper},  
%in: Algebra, Analysis and Beyond (Nowhere, 1973),   
%E.~Fox et al. (eds.),
%Lecture Notes in Math. 867, 
%Springer, Berlin, 1974, 115--124.

%\bibitem[N]{Nov} A. S. Novikov,
%\emph{Another fascinating article},  
%Uspekhi Mat. Nauk {23} (1980), no.~3, 112--134 (in Russian); 
%English transl.: Russian Math. Surveys 23 (1980), 572--595.




\end{thebibliography}
\end{document}